\documentclass[reqno]{amsart}

\usepackage{amssymb}
\usepackage{amsfonts}
\usepackage{txfonts}
\usepackage{dsfont}

\title{On pseudo-Riemannian Ricci-parallel  Lie groups which are not Einstein}
\author{Huihui An \and Zaili Yan$^{*}$}
\address[Huihui An]{ School of Mathematics, Liaoning Normal
University, Dalian, Liaoning Province, 116029, People's Republic of China}
\email[]{anhh@lnnu.edu.cn }
\address[Zaili Yan ]{School of Mathematics and Statistics, Ningbo University, Ningbo, Zhejiang Province, 315211,  People's Republic of China}
\email[]{yanzaili@nbu.edu.cn}
\thanks{$^{*}$Z. Yan is the corresponding
author and is supported by the Fundamental Research Funds for the Provincial Universities of Zhejiang and  K.C. Wong Magna Fund in Ningbo University.}
\date{}

\newtheorem{thm}{Theorem}[section]
\newtheorem{prop}[thm]{Proposition}
\newtheorem{lem}[thm]{Lemma}
\newtheorem{cor}[thm]{Corollary}
\theoremstyle{definition}

\newtheorem{example}[thm]{Example}

\newtheorem{rem}[thm]{Remark}

\begin{document}

\maketitle
\begin{abstract}
In this paper, we mainly study left invariant  pseudo-Riemannian Ricci-parallel metrics on connected Lie groups which are not Einstein. Following a result of Boubel and  B\'{e}rard Bergery, there are two typical  types of such metrics, which are  characterized by the minimal polynomial of the Ricci operator.  Namely, its form is either  $(X-\alpha)(X-\bar{\alpha})$ (type I), where $\alpha\in \mathbb{C}\setminus \mathbb{R}$, or $X^{2}$ (type II).
Firstly,  we obtain a complete description of Ricci-parallel metrics of type I. In particular, such a Ricci-parallel metric is uniquely determined by an Einstein metric and an invariant symmetric parallel complex structure up to isometry and scaling.  Then  we study Ricci-parallel  metric Lie algebras of  type II by using  double extension process. Surprisingly,  we find that every double extension of a metric Abelian Lie algebra is Ricci-parallel and the converse holds for Lorentz Ricci-parallel metric nilpotent Lie algebras of type II.
 Moreover,  we  construct  infinitely many new explicit examples of Ricci-parallel metric Lie algebras which are not Einstein.



\end{abstract}

\section{Introduction and main results}

The Ricci tensor $\mathrm{ric}$ of a (pseudo-)Riemannian manifold $(M,g)$ is an important quantity in differential  geometry and has been  extensively studied over the last one hundred  years.
$(M,g)$ is called Einstein if $\mathrm{ric}=cg$ holds for some constant $c\in \mathbb{R}$.
The study of homogeneous Riemannian Einstein metrics is particular interesting in this field and we refer the readers to \cite{ak75,bohm04,bk06,bk23,bl22,bl23,bwz04,dz79,heb98,lau10,Nikolayevsky11,nikonorov19,wz86}.
For the study of left invariant pseudo-Riemannian Einstein metrics on connected Lie groups,
see \cite{bt20,conti24,cbr21,cr20191,cr20193,cr2020,kath00,nomizu79,yd23,zy21}.
Here we collect some important results.
\begin{itemize}
    \item Several sufficient conditions for a compact homogeneous space admitting  a homogeneous   Riemannian Einstein metric are given in \cite{bohm04,bk23,bwz04,wz86}; every compact homogeneous space of dimension $\leq 11$ admits a homogeneous   Riemannian Einstein metric \cite{bk06}; every compact simple Lie group of dimension $\geq 4$ admits at least two non-isometric left invariant Riemannian Einstein metrics \cite{dz79}.
    \item Any connected homogeneous non-compact  Riemannian Einstein manifold with negative scalar curvature is isometric to a  Riemannian Einstein solvmanifold \cite{bl22,bl23}, i.e., a connected and simply connected solvable Lie group endowed with a left invariant  Riemannian Einstein metric, which is standard \cite{lau10}. Any solvable Lie group admits at most one left invariant  Riemannian Einstein metric up to isometry and scaling \cite{heb98}; every Riemannian Einstein solvmanifold determines a Ricci nilsoliton and conversely every Ricci nilsoliton can be extended to a rank-one Riemannian Einstein solvmanifold \cite{heb98,lau01,Nikolayevsky11}; a nilpotent Lie group admitting a left invariant Riemannian Einstein metric must be Abelian \cite{milnor76}.
    \item   Every non-compact simple Lie group of dimension $\geq 4$ admits at least two non-isometric left invariant pseudo-Riemannian Einstein metrics \cite{kath00}; every Riemannian Einstein solvmanifold admits a left invariant pseudo-Riemannian Einstein metric \cite{zy21}; some  pseudo-Riemannian Ricci nilsolitons can be extended to  pseudo-Riemannian Einstein solvmanifolds \cite{yd21}.
    \item  Every nilpotent Lie group of dimension $\leq 7$ has a left invariant  pseudo-Riemannian Ricci flat metric \cite{conti24,cbr21};  every nice nilpotent Lie group of dimension $\leq 11$ has a left invariant  pseudo-Riemannian Ricci flat metric \cite{conti24}; every nilpotent Lie group with large center  has a left invariant  pseudo-Riemannian Ricci flat metric \cite{yd23};  a nilpotent Lie  group endowed with a left invariant pseudo-Riemannian Einstein metric with non-zero scalar curvature is constructed in \cite{cr20191}.
\end{itemize}

Recall that a pseudo-Riemannian manifold $(M,g)$ (and the metric $g$) is called Ricci-parallel if $\nabla\mathrm{ric}=0$, where $\nabla$ denotes the Levi-Civita connection of $(M,g)$. Clearly, this notion is a natural generalization of Einstein manifolds, since  every Einstein manifold is Ricci-parallel.
 Gray in \cite{gray78} proved that a Riemannnian Ricci-parallel manifold  is locally a Riemannian product of Einstein manifolds. On the other hand, in the pseudo-Riemannian case, Boubel and B\'{e}rard  Bergery \cite{bb01} showed that there are two types of  pseudo-Riemannian non-Einstein Ricci-parallel manifolds.
More precisely, let $(M,g)$ be an indecomposable  pseudo-Riemannian non-Einstein Ricci-parallel manifold, then its Ricci operator $\mathrm{Ric}$ must satisfy one of the following two conditions (see the Main Theorem of \cite{bb01}).
\begin{itemize}
  \item  (Type I) The minimal polynomial of $\mathrm{Ric}$ is of the form $(X-\alpha)(X-\bar{\alpha})$, where $\alpha=\lambda+\mathbf{i}\mu$, $\lambda,\mu\in \mathbb{R}$, $\mu\neq 0$.
  \item (Type II)  The minimal polynomial of $\mathrm{Ric}$ has the form $X^{2}$, i.e., $\mathrm{Ric}^{2}=0$ but $\mathrm{Ric}\neq 0$.
\end{itemize}

Furthermore, Boubel and B\'{e}rard Bergery \cite{bb01} proved that a family of pseudo-Riemannian symmetric spaces due to Cahen and  Parker \cite{cp70} belongs to type II, and every pseudo-Riemannian manifold $(M,g)$ of type I admits an Einstein metric whose Levi-Civita connection as the same as $g$.
Up to now, to the best knowledge of the authors,  pseudo-Riemannian manifolds of these two types are less known and a classification of such manifolds seems to be unreachable, so it would be interesting and important to find more examples.

This paper is devoted to the study of left invariant pseudo-Riemannian Ricci-parallel metrics of types I and II on connected Lie groups.
Assume $(\mathfrak{g},\langle\cdot,\cdot\rangle)$ is a metric Lie algebra, which stands for a connected Lie group $G$ with Lie algebra $\mathfrak{g}$ endowed with a left invariant pseudo-Riemannian metric $g$ generated by $\langle\cdot,\cdot\rangle$.
For type I,  inspired by Proposition 1 and Proposition 3 of \cite{bb01}, we obtain a complete description of such Ricci-parallel metric Lie algebras.
\begin{thm}\label{1-thm-main1}
Given a complex number $\alpha=\lambda+\mathbf{i}\mu\in \mathbb{C}\setminus \mathbb{R}$, $\lambda, \mu\in \mathbb{R}$,
and let $(\mathfrak{g},\langle\cdot,\cdot\rangle)$ be a Ricci-parallel metric Lie algebra.
Then the minimal polynomial of $\mathrm{Ric}$ is of the form $(X-\alpha)(X-\bar{\alpha})$ if and only if
there exists a unique pair $(\langle\cdot,\cdot\rangle',\mathbf{J})$ on $\mathfrak{g}$ determined by the following equation
\begin{equation}\label{main1-equ}
\langle\cdot,\cdot\rangle=\frac{1}{\lambda^{2}+\mu^{2}}
\left[\lambda\langle\cdot,\cdot\rangle'-\mu\langle\cdot,\mathbf{J}(\cdot)\rangle'\right],
\end{equation}
where $\langle\cdot,\cdot\rangle'$ is an invariant Einstein metric  on $\mathfrak{g}$ with Einstein constant 1 and $\mathbf{J}:\mathfrak{g}\rightarrow \mathfrak{g}$  is an invariant symmetric parallel (with respect to $\langle\cdot,\cdot\rangle'$) complex structure on $\mathfrak{g}$.
\end{thm}

For type II, we provide a sufficient condition to obtain such  Ricci-parallel metric Lie algebras. This method is called double extension, which was first  introduced by   Medina and   Revoy \cite{mr85} in the context of ad-invariant metrics.
 In practice, a  metric Lie algebra $(\mathfrak{g},\langle\cdot,\cdot\rangle)$  is said to be a double extension of another metric Lie algebra
$(\mathfrak{g}_{0},\langle\cdot,\cdot\rangle_{0})$ if it has the following form:
$$\mathfrak{g}=\mathbb{R}u+\mathbb{R}v+\mathfrak{g}_{0},$$
 $$[u,e]=D(e)+\langle L,e\rangle_{0}v,\quad [e,e']=[e,e']_{0}+\langle K(e),e'\rangle_{0}v,
\quad [\mathfrak{g},v]=0,\quad e,e'\in \mathfrak{g}_{0},$$
$$\langle u,v\rangle=1,\quad \langle u,u\rangle=\langle v,v\rangle=0,$$
$$\langle u,\mathfrak{g}_{0}\rangle=\langle v,\mathfrak{g}_{0}\rangle=0, \quad
\langle \cdot,\cdot\rangle|_{\mathfrak{g}_{0}}=\langle\cdot,\cdot\rangle_{0},$$
where $D,K\in \mathrm{End}(\mathfrak{g}_{0})$, $L\in \mathfrak{g}_{0}$,
$[\cdot,\cdot]$ and $[\cdot,\cdot]_{0}$ denote the Lie brackets of $\mathfrak{g}$ and $\mathfrak{g}_{0}$, respectively.

Now let $\Delta\in \mathfrak{g}_{0}$ be a vector given by
\begin{equation*}
  \langle \Delta, e\rangle_{0}=-\frac{1}{2}\mathrm{tr}((D+D^{*})\circ \mathrm{ad}_{0}\,e)
  +\frac{1}{2}\langle (D-K)(Z_{0}),e\rangle_{0}
 -\frac{1}{4}\mathrm{tr}(KS^{0}_{e}),\quad \forall e\in \mathfrak{g}_{0},
\end{equation*}
where $Z_{0}$ denotes the mean curvature vector of $(\mathfrak{g}_{0},\langle\cdot,\cdot\rangle_{0})$,
i.e., $\langle Z_{0},f\rangle_{0}=\mathrm{tr}(\mathrm{ad}_{0}\,f)$, $\forall f\in \mathfrak{g}_{0}$,
and $S^{0}_{e}:\mathfrak{g}_{0}\rightarrow \mathfrak{g}_{0}$ is defined by $S^{0}_{e}(f)=(\mathrm{ad}_{0}\,f)^{*}(e)$ for all $e\in \mathfrak{g}_{0}$.
Here $\mathrm{ad}_{0}$ denotes the adjoint map of $\mathfrak{g}_{0}$
and $A^{*}:\mathfrak{g}_{0}\rightarrow \mathfrak{g}_{0}$ denotes the  adjoint map of $A:\mathfrak{g}_{0}\rightarrow \mathfrak{g}_{0}$ with  respect to $\langle\cdot,\cdot\rangle_{0}$. Moreover, let
$$\Gamma=-\frac{1}{2}\mathrm{tr}(D^{2})-\frac{1}{2}\mathrm{tr}(D^{*}D)
  -\frac{1}{4}\mathrm{tr}(K^{2})+\langle L,Z_{0}\rangle,$$
we have
\begin{thm}\label{1-thm-main2}
Let $(\mathfrak{g},\langle\cdot,\cdot\rangle)$ be a double extension of $(\mathfrak{g}_{0},\langle\cdot,\cdot\rangle_{0})$ as described  above.
Assume that $(\mathfrak{g}_{0},\langle\cdot,\cdot\rangle_{0})$ is Ricci flat and $\Delta=0$,  then $(\mathfrak{g},\langle\cdot,\cdot\rangle)$ is Ricci-parallel.
Moreover, if in addition that $\Gamma\neq0$, then $(\mathfrak{g},\langle\cdot,\cdot\rangle)$ is Ricci-parallel with $\mathrm{Ric}\neq0$ and $\mathrm{Ric}^{2}=0$.
\end{thm}
Notice that $\Delta=0$ holds for every metric Abelian Lie algebra, so we have
\begin{cor}
Every double extension of a metric Abelian Lie algebra is Ricci-parallel.
\end{cor}
Conversely, we prove the following result.
\begin{thm}\label{1-thm-main3}
Every Lorentz Ricci-parallel  metric nilpotent Lie algebra of type  II
is a double extension of Euclidean metric Abelian Lie algebra.
\end{thm}

 Finally, we present an interesting class of Ricci-parallel metric Lie algebras, which is inspired by the  structure of  even-dimensional ad-invariant metric  Lie algebras \cite{bordemann97}.
\begin{example}\label{1-exam}
Let $\mathfrak{D}$ be  an $n$-dimensional  Lie algebra with the dual vector space $\mathfrak{D}^{*}$,
and $\theta\in Z^{2}(\mathfrak{D},\mathfrak{D}^{*})$ be a cocycle.
Let $\mathfrak{g}=\mathfrak{D}+\mathfrak{D}^{*}$ be the central extension of $\mathfrak{D}$ by $\theta$, namely, it is  defined  by
$$[x_{1}+f_{1},x_{2}+f_{2}]=[x_{1},x_{2}]_{\mathfrak{D}}+\theta(x_{1},x_{2}), \quad x_{1}, x_{2}\in \mathfrak{D}, f_{1},f_{2}\in \mathfrak{D}^{*}.$$
We claim that the natural  metric $\langle\cdot,\cdot\rangle$ of signature $(n,n)$ on $\mathfrak{g}=\mathfrak{D}+\mathfrak{D}^{*}$ defined by
$$\langle x_{1}+f_{1}, x_{2}+f_{2}\rangle=f_{1}(x_{2})+f_{2}(x_{1})$$
is Ricci-parallel. In particular, it is Ricci flat if and only if its Killing form vanishes.
\end{example}

This paper is organized as follows.
In Section \ref{sec2}, we recall some basic facts about metric Lie algebras.
In Section \ref{sec3}, we  study the Ricci-parallel metric Lie algebras of type I and prove Theorem \ref{1-thm-main1}.
In Section \ref{sec4}, we study the  Ricci-parallel metric Lie algebras of type II and prove Theorem \ref{1-thm-main2}.
In Section \ref{sec5}, we study Lorentz Ricci-parallel  metric Lie algebras of type II and prove Theorem \ref{1-thm-main3}.

\medskip
\section{Metric Lie algebras}\label{sec2}
In this section, we recall some basic facts about  left invariant pseudo-Riemannian metrics on connected Lie groups.
Let $G$ be a connected $n$-dimensional Lie group and $\mathfrak{g}$ be the  Lie algebra of $G$ consisting of left invariant vector fields on  $G$.  Then each left invariant pseudo-Riemannian metric  $\langle\cdot,\cdot\rangle$ on $G$ is in fact a non-degenerate symmetric
bilinear form $\langle\cdot,\cdot\rangle$ on $\mathfrak{g}$.
For convenience, we call such a pair $(\mathfrak{g}, \langle\cdot,\cdot\rangle)$ a metric Lie algebra.
 The metric $\langle\cdot,\cdot\rangle$ is said to be of signature $(p,q)$  if there exists a basis $\{e_{1},e_{2},\ldots,e_{n}\}$ of $\mathfrak{g}$,  such that the metric matrix $(\langle e_{i},e_{j}\rangle)=\mathrm{diag} (\underbrace{-1,\cdots,-1}_{p},\underbrace{1,\cdots,1}_{q})$, $p+q=n$.
In particular, it is called Euclidean (resp. Lorentz) if the signature is $(0,n)$ (resp. $(1,n-1)$).

Now let $\nabla$ be the Levi-Civita connection of $(\mathfrak{g},\langle\cdot,\cdot\rangle)$.
Then for   any $x, y\in \mathfrak{g}$, the Koszul formula gives
\begin{eqnarray*}
  \nabla_{x}y &=&\frac{1}{2}\left([x,y]-(\mathrm{ad}\,x)^{*}y-(\mathrm{ad}\,y)^{*}x\right)\\
  &=& \frac{1}{2}[x,y]+U(x,y),
\end{eqnarray*}
where $(\mathrm{ad}\,x)^{*}$ denotes the adjoint map of $\mathrm{ad}\,x$ with respect to $\langle\cdot,\cdot\rangle$,
and $U: \mathfrak{g}\times \mathfrak{g}\rightarrow \mathfrak{g}$ is defined by
\begin{equation*}
 \langle U(x,y),z\rangle=\frac{1}{2}\left(\langle[z,x],y\rangle+\langle x,[z,y]\rangle\right),
 \quad \forall x,y,z\in \mathfrak{g}.
\end{equation*}
The Riemannian curvature tensor $R$ is defined in terms of $\nabla$ by the following form:
\begin{equation}\label{}
  R(x,y)z=[\nabla_{x},\nabla_{y}]z-\nabla_{[x,y]}z=\nabla_{x}\nabla_{y}z-\nabla_{y}\nabla_{x}z-\nabla_{[x,y]}z,\quad x,y,z\in \mathfrak{g}. \nonumber
\end{equation}
The Ricci tensor $\mathrm{ric}$ and the Ricci  operator $\mathrm{Ric}$ are defined by
\begin{equation}\label{}
  \mathrm{ric}(x,y)=\mathrm{trace} (\xi\longrightarrow R(\xi,x)y),\quad
  \langle\mathrm{Ric}(x),y\rangle=\mathrm{ric}(x,y),
  \quad \forall x,y\in \mathfrak{g}.  \nonumber
\end{equation}
To give the formula of Ricci tensor $\mathrm{ric}$, we  assume $\{e_{i}\}$ is an orthonormal basis of $(\mathfrak{g},\langle\cdot,\cdot\rangle)$,
that is, $\langle e_{i},e_{i}\rangle=\varepsilon_{i}\in \{1,-1\}$, $\langle e_{i},e_{j}\rangle=0$ for  $i\neq j$.
Then  one  easily deduces that (see also \cite{oneill1983book})
  \begin{eqnarray}\label{2-formula-ric}
 &&\mathrm{ric}(x,y)=\sum_{i=1}^{n}\langle R(e_{i},x)y,e_{i}\rangle \varepsilon_{i} \nonumber\\
 &=& -\frac{1}{2}\mathbf{K}(x,y)-\frac{1}{2}(\langle [Z_{\mathfrak{g}},x],y\rangle+\langle [Z_{\mathfrak{g}},y],x\rangle) \nonumber\\
  & &-\frac{1}{2}\sum^{n}_{i=1}\langle[x,e_{i}],[y,e_{i}]\rangle\varepsilon_{i}+ \frac{1}{4}\sum^{n}_{i,j=1}\langle[e_{i},e_{j}],x\rangle\langle[e_{i},e_{j}],y\rangle\varepsilon_{i}\varepsilon_{j} \nonumber\\
  &=&-\frac{1}{2}\mathbf{K}(x,y)-\frac{1}{2}(\langle [Z_{\mathfrak{g}},x],y\rangle+\langle [Z_{\mathfrak{g}},y],x\rangle)-\frac{1}{2}\mathrm{tr}((\mathrm{ad}\,x)^{*}\mathrm{ad}\,y)
  +\frac{1}{4}\mathrm{tr}(S^{*}_{x}S_{y}),
 \end{eqnarray}
$\forall x,y\in \mathfrak{g}$, where $\mathbf{K}$ is the Killing form of $\mathfrak{g}$,
$Z_{\mathfrak{g}}$ is the mean curvature vector in $\mathfrak{g}$ satisfying
$\langle Z_{\mathfrak{g}},x\rangle=\mathrm{tr}(\mathrm{ad}\,x)$ for all $x\in \mathfrak{g}$,
the operator $S_{x}:\mathfrak{g}\rightarrow \mathfrak{g}$ is defined by $S_{x}(y)=(\mathrm{ad}\,y)^{*}(x)$ for all $x\in \mathfrak{g}$.
 Notice that $Z_{\mathfrak{g}}=0$ if and only if $\mathfrak{g}$ is unimodular,
if and only if $\mathrm{tr}(\mathrm{ad}\, x)=0$ holds for all $x\in \mathfrak{g}$.
Moreover, one can mention that the center $C(\mathfrak{g})\subset \ker(S_{x})$ for all $x\in \mathfrak{g}$.

We say that $(\mathfrak{g},\langle\cdot,\cdot\rangle)$ is Einstein if there exists $c\in \mathbb{R}$ such that $\mathrm{ric}=c\langle\cdot,\cdot\rangle$. In particular, $(\mathfrak{g},\langle\cdot,\cdot\rangle)$ is Ricci flat if $\mathrm{ric}$ is identically zero.
Moreover, $(\mathfrak{g},\langle\cdot,\cdot\rangle)$ is called Ricci-parallel if $\nabla \mathrm{ric}=0$.
The following result is basic but useful.
\begin{prop}\label{2-prop-parallel}
$(\mathfrak{g},\langle\cdot,\cdot\rangle)$ is Ricci-parallel if and only if
$\mathrm{Ric}\circ\nabla_{x}=\nabla_{x}\circ\mathrm{Ric}$ holds for all $x\in \mathfrak{g}$, where $\mathrm{Ric}$ denotes the Ricci operator of $(\mathfrak{g},\langle\cdot,\cdot\rangle)$.
\end{prop}
\begin{proof}
Notice that $\nabla \mathrm{ric}=0$ if and only if for all $x,y,z\in \mathfrak{g}$,
\begin{eqnarray*}
0&=&(\nabla_{x} \mathrm{ric})(y,z)\\
 &=& -\mathrm{ric}(\nabla_{x}y,z)-\mathrm{ric}(y,\nabla_{x}z)\\
 &=&-\langle(\mathrm{Ric}\circ \nabla_{x})(y),z\rangle-\langle \mathrm{Ric}(y),\nabla_{x}z\rangle\\
 &=&-\langle(\mathrm{Ric}\circ \nabla_{x})(y),z\rangle+\langle(\nabla_{x}\circ\mathrm{Ric})(y),z\rangle\\
 &=&\left\langle(\nabla_{x}\circ\mathrm{Ric}-\mathrm{Ric}\circ \nabla_{x})(y),z\right\rangle.
 \end{eqnarray*}
 Hence $\nabla \mathrm{ric}=0$ is equivalent to saying that $\mathrm{Ric}\circ\nabla_{x}=\nabla_{x}\circ\mathrm{Ric}$, $\forall x\in \mathfrak{g}$.
\end{proof}

Now we are in a position to prove the claim of Example \ref{1-exam}.

\textbf{Proof of the claim of Example \ref{1-exam}}.
By definitions and the fact that $\mathfrak{D}^{*}\subset C(\mathfrak{g})$, the center of $\mathfrak{g}$,
 we easily obtain that
\begin{equation*}
 Z_{\mathfrak{g}}\in \mathfrak{D}^{*}, \quad
  (\mathrm{ad}\,x)^{*}(\mathfrak{g})\subset \mathfrak{D}^{*},\quad
  S_{x}(\mathfrak{g})\subset \mathfrak{D}^{*}, \quad
  S^{*}_{x}(\mathfrak{D}^{*})=0,
  \quad U(x,y)\in \mathfrak{D}^{*},\quad \forall x,y\in \mathfrak{g}.
\end{equation*}
So by formula \eqref{2-formula-ric}, we have
\begin{equation*}
\mathrm{ric}(x,y)=-\frac{1}{2}\mathbf{K}(x,y), \forall x,y\in \mathfrak{g}.
\end{equation*}
Since $\mathbf{K}(\mathfrak{g},\mathfrak{D}^{*})=0$,
so for all $x,y,z\in \mathfrak{g}$, we obtain
\begin{eqnarray*}
(\nabla_{x} \mathrm{ric})(y,z)&=& -\mathrm{ric}(\nabla_{x}y,z)-\mathrm{ric}(y,\nabla_{x}z)\\
 &=&\frac{1}{2}\mathbf{K}(\nabla_{x}y,z)+\frac{1}{2}\mathbf{K}(y,\nabla_{x}z)\\
 &=&\frac{1}{4}\mathbf{K}([x,y],z)+\frac{1}{4}\mathbf{K}(y,[x,z])\\
 &=&0.
 \end{eqnarray*}
 Hence $(\mathfrak{g},\langle\cdot,\cdot\rangle)$ is Ricci-parallel. In particular, it is Ricci flat if and only if its Killing form $\mathbf{K}=0$.\qed

At the  last of  this section, we briefly recall the notion of  ad-invariant metrics on
Lie algebras. For more results, we refer the readers to an excellent survey \cite{ovando16} and the
references therein.
A  metric $\langle\cdot,\cdot\rangle$ on a Lie algebra $\mathfrak{g}$ is called ad-invariant if it satisfies
\begin{equation*}
 \langle [x,y],z\rangle+\langle y,[x,z]\rangle=0
\end{equation*}
for all $x,y,z\in \mathfrak{g}$.
By definition, the Levi-Civita connection $\nabla$, Riemannian curvature tensor $R$ and Ricci tensor $\mathrm{ric}$ of  $(\mathfrak{g}, \langle\cdot,\cdot\rangle)$ are given by
\begin{equation*}
\nabla_{x}y=\frac{1}{2}[x,y], \quad R(x,y)z=-\frac{1}{4}[[x,y],z], \quad  \mathrm{ric}(x,y)=-\frac{1}{4}\mathbf{K}(x,y), \quad \forall x,y\in \mathfrak{g},
\end{equation*}
respectively.
In particular, $(\mathfrak{g}, \langle\cdot,\cdot\rangle)$ is Ricci flat if $\mathfrak{g}$ is nilpotent.
Clearly, $\nabla R=\nabla \mathrm{ric}=0$. So we have
\begin{prop}
Every ad-invariant metric Lie algebra $(\mathfrak{g}, \langle\cdot,\cdot\rangle)$ is Ricci-parallel.
\end{prop}
\begin{rem}
 In fact, the pseudo-Riemannian manifold $(G, \langle\cdot,\cdot\rangle)$ corresponding to an ad-invariant metric Lie algebra $(\mathfrak{g}, \langle\cdot,\cdot\rangle)$ is a pseudo-Riemannian symmetric space.
\end{rem}

\section{Ricci-parallel metric Lie algebras of type I}\label{sec3}

In this section, we study the Ricci-parallel metric Lie algebras of type I.  First, we state the proof of Theorem \ref{1-thm-main1}.

\textbf{Proof of Theorem \ref{1-thm-main1}}.
Let $(\mathfrak{g},\langle\cdot,\cdot\rangle)$ be a Ricci-parallel metric Lie algebra of type I.

$"\Rightarrow"$: Assume the minimal polynomial of $\mathrm{Ric}$ is of the form: $(X-\alpha)(X-\bar{\alpha})$, where $\alpha=\lambda+\mathbf{i}\mu$, $\lambda,\mu\in \mathbb{R}$, $\mu\neq 0$.
Then let $\mathbf{J}=\frac{1}{\mu}(\mathrm{Ric}-\lambda \mathrm{Id})$ and
$$\langle\cdot,\cdot\rangle'=\mathrm{Re}(\alpha h)=\lambda\langle\cdot,\cdot\rangle+\mu\langle\cdot,\mathbf{J}(\cdot)\rangle,$$
where $h=\langle\cdot,\cdot\rangle-\mathbf{i}\langle\cdot,\mathbf{J}(\cdot)\rangle$.
According to the proof of Proposition 1 of \cite{bb01},
the Levi-Civita connections of $h$ and $\langle\cdot,\cdot\rangle'$ are the same as $\langle\cdot,\cdot\rangle$.
Moreover, $\langle\cdot,\cdot\rangle'$ is Einstein with Einstein constant 1 and $\mathbf{J}$
is  a symmetric parallel  complex structure  with respect to $\langle\cdot,\cdot\rangle'$.
Notice that
$$\langle\cdot,\mathbf{J}(\cdot)\rangle'=
\lambda\langle\cdot,\mathbf{J}(\cdot)\rangle-\mu\langle\cdot,\cdot\rangle,$$
we obtain
$$\langle\cdot,\cdot\rangle=\frac{1}{\lambda^{2}+\mu^{2}}
\left[\lambda\langle\cdot,\cdot\rangle'-\mu\langle\cdot,\mathbf{J}(\cdot)\rangle'\right].$$

$"\Leftarrow"$:
Suppose that $(\mathfrak{g},\langle\cdot,\cdot\rangle')$ is an Einstein  metric Lie algebra with Einstein constant 1 and $\mathbf{J}$
is  a symmetric parallel  complex structure on $\mathfrak{g}$ with respect to $\langle\cdot,\cdot\rangle'$.
By definition, it is easily seen that $\langle\cdot,\cdot\rangle$ induces the same Levi-Civita connection as $\langle\cdot,\cdot\rangle'$, so the same Ricci tensor $\mathrm{ric}=\langle\cdot,\cdot\rangle'$. This also implies the uniqueness of $\langle\cdot,\cdot\rangle'$ and $\mathbf{J}$.
Consequently, the Ricci operator $\mathrm{Ric}$ of $(\mathfrak{g},\langle\cdot,\cdot\rangle)$  satisfies
$\langle \mathrm{Ric}(\cdot),\cdot\rangle=\mathrm{ric}=\langle\cdot,\cdot\rangle'$.
On the other hand, since
$$\langle\cdot,\cdot\rangle=\frac{1}{\lambda^{2}+\mu^{2}}
\left[\lambda\langle\cdot,\cdot\rangle'-\mu\langle\cdot,\mathbf{J}(\cdot)\rangle'\right]$$
and then
$$\langle\cdot,\mathbf{J}(\cdot)\rangle=\frac{1}{\lambda^{2}+\mu^{2}}
\left[\lambda\langle\cdot,\mathbf{J}(\cdot)\rangle'+\mu\langle\cdot,\cdot\rangle'\right],$$
we obtain
$$\langle\cdot,\cdot\rangle'=
\lambda\langle\cdot,\cdot\rangle+\mu\langle\cdot,\mathbf{J}(\cdot)\rangle.$$
Therefore $\mathrm{Ric}=\lambda\mathrm{Id}+\mu \mathbf{J}$, which completes the proof of the theorem.\qed

We say two metric Lie algebras $(\mathfrak{g}_{1},\langle\cdot, \cdot\rangle_{1})$ and $(\mathfrak{g}_{2},\langle\cdot, \cdot\rangle_{2})$ are isometric if there exists a Lie algebra isomorphism $\varphi: \mathfrak{g}_{1}\rightarrow \mathfrak{g}_{2}$ such that
$\langle\cdot, \cdot\rangle _{1}=\langle\varphi(\cdot), \varphi(\cdot)\rangle_{2}.$
The following result is a direct corollary of Theorem \ref{1-thm-main1}.
\begin{cor}
Given a complex number $\alpha=\lambda+\mathbf{i}\mu\in \mathbb{C}\setminus \mathbb{R}$,
 $\lambda,\mu\in \mathbb{R}$, two Ricci-parallel metric Lie algebras $(\mathfrak{g}_{1},\langle\cdot, \cdot\rangle_{1})$ and $(\mathfrak{g}_{2},\langle\cdot, \cdot\rangle_{2})$ with the same minimal polynomial $(X-\alpha)(X-\bar{\alpha})$ are isometric if and only if
the related triples $(\mathfrak{g}_{1},\langle\cdot, \cdot\rangle'_{1}, \mathbf{J}_{1})$ and $(\mathfrak{g}_{2},\langle\cdot, \cdot\rangle'_{2},\mathbf{J}_{2})$ determined by \eqref{main1-equ} in Theorem \ref{1-thm-main1} are isometric.
Namely, there exists a Lie algebra isomorphism $\varphi: \mathfrak{g}_{1}\rightarrow \mathfrak{g}_{2}$ such that
$\langle\cdot, \cdot\rangle'_{1}=\langle\varphi(\cdot), \varphi(\cdot)\rangle'_{2}$ and $\mathbf{J}_{2}=\varphi\circ\mathbf{J}_{1}\circ \varphi^{-1}$.
\end{cor}

Next, we provide a method to construct Ricci-parallel metric Lie algebras of type I.
\begin{prop}\label{3-prop-complex}
Suppose that $(\mathfrak{g}_{0},\langle\cdot,\cdot\rangle_{0})$ is an Einstein metric Lie algebra with Einstein constant $c$. Let $\mathfrak{g}=\mathfrak{g}_{0}+\mathbf{i}\mathfrak{g}_{0}$ be the complexification of $\mathfrak{g}_{0}$ and $\mathbf{J}$ be the natural almost complex structure on $\mathfrak{g}$ defined by
$$\mathbf{J}(x+\mathbf{i}y)=-y+\mathbf{i}x,\quad \forall x,y\in \mathfrak{g}_{0}.$$
Regarding $\mathfrak{g}$ as a real Lie algebra and  define a metric $\langle\cdot,\cdot\rangle'$ on $\mathfrak{g}$ as follows:
$$\langle x_{1}+\mathbf{i}y_{1},x_{2}+\mathbf{i}y_{2}\rangle'
=\langle x_{1},x_{2}\rangle_{0}-\langle y_{1},y_{2}\rangle_{0},
\quad \forall x_{1},x_{2}, y_{1},y_{2}\in \mathfrak{g}_{0}.$$
Then $(\mathfrak{g},\langle\cdot,\cdot\rangle')$ is an Einstein metric Lie algebra with Einstein constant $2c$ and $\mathbf{J}$  is an symmetric parallel complex structure on $(\mathfrak{g},\langle\cdot,\cdot\rangle')$.
Consequently, if $c\neq 0$, then the metric $\langle\cdot,\cdot\rangle$ on $\mathfrak{g}$ defined by
\begin{eqnarray*}
  \langle x_{1}+\mathbf{i}y_{1},x_{2}+\mathbf{i}y_{2}\rangle &=&
   \frac{2c}{\lambda^{2}+\mu^{2}}[\lambda\langle x_{1}+\mathbf{i}y_{1},x_{2}+\mathbf{i}y_{2}\rangle'
   -\mu\langle x_{1}+\mathbf{i}y_{1},\mathbf{J}(x_{2}+\mathbf{i}y_{2})\rangle'] \\
   &=&\frac{2c}{\lambda^{2}+\mu^{2}}[\lambda\langle x_{1},x_{2}\rangle_{0}-\lambda\langle y_{1},y_{2}\rangle_{0}+\mu\langle x_{1},y_{2}\rangle_{0}+\mu\langle y_{1},x_{2}\rangle_{0}],
\end{eqnarray*}
$\forall x_{1},x_{2}, y_{1},y_{2}\in \mathfrak{g}_{0}$, is Ricci-parallel, and the minimal polynomial of its Ricci operator is of the form $(X-\lambda-\mathbf{i}\mu)(X-\lambda+\mathbf{i}\mu)$, where $\lambda,\mu\in \mathbb{R}$ and $\mu\neq 0$.
\end{prop}
\begin{proof}
Let $\nabla$ $(\widetilde{\nabla})$, $R$ $(\widetilde{R})$ and $\mathrm{ric}$ $(\widetilde{\mathrm{ric}})$ be the Levi-Civita connection, Riemannian curvature tensor and Ricci tensor of $(\mathfrak{g}_{0},\langle\cdot,\cdot\rangle_{0})$ $((\mathfrak{g},\langle\cdot,\cdot\rangle'))$, respectively.
By a direct calculation, we obtain
$$\widetilde{\nabla}_{x_{1}+\mathbf{i}y_{1}}(x_{2}+\mathbf{i}y_{2})
=\nabla_{x_{1}}x_{2}-\nabla_{y_{1}}y_{2}+\mathbf{i}(\nabla_{x_{1}}y_{2}+\nabla_{y_{1}}x_{2}),
\quad \forall x_{1},x_{2}, y_{1},y_{2}\in \mathfrak{g}_{0}.$$
Hence
\begin{gather*}
  \widetilde{R}(x,y)z=R(x,y)z, \quad \widetilde{R}(x,\mathbf{i}y)z=\mathbf{i}R(x,y)z,\quad \widetilde{R}(\mathbf{i}x,\mathbf{i}y)z=-R(x,y)z,\\
  \widetilde{R}(x,y)(\mathbf{i}z)=\mathbf{i}R(x,y)z, \quad \widetilde{R}(x,\mathbf{i}y)(\mathbf{i}z)=-R(x,y)z,\quad \widetilde{R}(\mathbf{i}x,\mathbf{i}y)(\mathbf{i}z)=-\mathbf{i}R(x,y)z
\end{gather*}
hold for all $x,y,z\in \mathfrak{g}_{0}$.
Therefore, for all $x,y\in \mathfrak{g}_{0}$,
\begin{gather*}
 \widetilde{\mathrm{ric}}(x,y)=2\mathrm{ric}(x,y)=2c\langle x,y\rangle_{0}=2c\langle x,y\rangle',  \\
  \widetilde{\mathrm{ric}}(x,\mathbf{i}y)=0, \\
  \widetilde{\mathrm{ric}}(\mathbf{i}x,\mathbf{i}y)=-2\mathrm{ric}(x,y)=-2c\langle x,y\rangle_{0}=2c\langle \mathbf{i}x,\mathbf{i}y\rangle'.
\end{gather*}
This implies that $(\mathfrak{g},\langle\cdot,\cdot\rangle')$ is an Einstein metric Lie algebra with Einstein constant $2c$.

Furthermore, notice that for all $x_{1},x_{2}, y_{1},y_{2}\in \mathfrak{g}_{0}$,
\begin{eqnarray*}
\mathbf{J}\left(\widetilde{\nabla}_{x_{1}+\mathbf{i}y_{1}}(x_{2}+\mathbf{i}y_{2})\right) &=& \mathbf{J}\left(\nabla_{x_{1}}x_{2}-\nabla_{y_{1}}y_{2}+\mathbf{i}(\nabla_{x_{1}}y_{2}+\nabla_{y_{1}}x_{2})\right) \\
&=& -(\nabla_{x_{1}}y_{2}+\nabla_{y_{1}}x_{2})+\mathbf{i}(\nabla_{x_{1}}x_{2}-\nabla_{y_{1}}y_{2})
\end{eqnarray*}
and
\begin{eqnarray*}
   \widetilde{\nabla}_{x_{1}+\mathbf{i}y_{1}}\mathbf{J}(x_{2}+\mathbf{i}y_{2}) &=&
    \widetilde{\nabla}_{x_{1}+\mathbf{i}y_{1}}(-y_{2}+\mathbf{i}x_{2})  \\
   &=& -(\nabla_{y_{1}}x_{2}+\nabla_{x_{1}}y_{2})+\mathbf{i}(\nabla_{x_{1}}x_{2}-\nabla_{y_{1}}y_{2}).
\end{eqnarray*}
We obtain that
$\mathbf{J}\circ\widetilde{\nabla}_{z}=\widetilde{\nabla}_{z}\circ\mathbf{J}$ holds for all $z\in \mathfrak{g}$.
 So by Proposition \ref{2-prop-parallel}, $\mathbf{J}$ is parallel and consequently a complex structure on $(\mathfrak{g},\langle\cdot,\cdot\rangle')$ (see \cite{nn57}).
 Clearly, $\mathbf{J}$ is symmetric with respect to $\langle\cdot,\cdot\rangle'$.
 The final part is a direct corollary of Theorem \ref{1-thm-main1}.
\end{proof}
\begin{rem}
We don't known whether the converse of Proposition \ref{3-prop-complex} is true.
\end{rem}

Moreover, we have the following result.
\begin{thm}
Suppose that $(\mathfrak{g}_{0},\langle\cdot,\cdot\rangle_{0})$ is a Ricci-parallel metric Lie algebra
and $\mathfrak{g}=\mathfrak{g}_{0}+\mathbf{i}\mathfrak{g}_{0}$ is  the complexification of $\mathfrak{g}_{0}$.
Regarding $\mathfrak{g}$ as a real Lie algebra and  then the metric $\langle\cdot,\cdot\rangle'$ defined on $\mathfrak{g}$ as in Proposition \ref{3-prop-complex} by
$$\langle x_{1}+\mathbf{i}y_{1},x_{2}+\mathbf{i}y_{2}\rangle'
=\langle x_{1},x_{2}\rangle_{0}-\langle y_{1},y_{2}\rangle_{0},
\quad \forall x_{1},x_{2}, y_{1},y_{2}\in \mathfrak{g}_{0},$$
is also Ricci-parallel.
\end{thm}
\begin{proof}
Keep notation as in Proposition \ref{3-prop-complex}.
Let $\mathrm{Ric}$ and $\widetilde{\mathrm{Ric}}$ be  Ricci operators  of $(\mathfrak{g}_{0},\langle\cdot,\cdot\rangle_{0})$ and $(\mathfrak{g},\langle\cdot,\cdot\rangle')$, respectively.
By the arguments in the proof of Proposition \ref{3-prop-complex}, we obtain that
$\widetilde{\mathrm{Ric}}(x+\mathbf{i}y)=2\mathrm{Ric}(x)+2\mathbf{i}\mathrm{Ric}(y)$, $\forall x,y\in \mathfrak{g}_{0}$. A direct calculation shows the following equalities:
\begin{eqnarray*}
   && \left(\widetilde{\nabla}_{x_{1}+\mathbf{i}y_{1}}\circ \widetilde{\mathrm{Ric}}\right)(x_{2}+\mathbf{i}y_{2}) \\
   &=& \widetilde{\nabla}_{x_{1}+\mathbf{i}y_{1}}(2\mathrm{Ric}(x_{2})+2\mathbf{i}\mathrm{Ric}(y_{2})) \\
   &=&  2\left(\nabla_{x_{1}}\mathrm{Ric}(x_{2})-\nabla_{y_{1}}\mathrm{Ric}(y_{2})
   +\mathbf{i}\left(\nabla_{x_{1}}\mathrm{Ric}(y_{2})+\nabla_{y_{1}}\mathrm{Ric}(x_{2})\right)\right),
\end{eqnarray*}
and
\begin{eqnarray*}
&& \left(\widetilde{\mathrm{Ric}}\circ\widetilde{\nabla}_{x_{1}+\mathbf{i}y_{1}}\right)(x_{2}+\mathbf{i}y_{2}) \\
&=& \widetilde{\mathrm{Ric}}\left(\nabla_{x_{1}}x_{2}-\nabla_{y_{1}}y_{2}+\mathbf{i}(\nabla_{x_{1}}y_{2}+\nabla_{y_{1}}x_{2})\right) \\
&=&  2\mathrm{Ric}\left(\nabla_{x_{1}}x_{2}-\nabla_{y_{1}}y_{2}\right)
+2\mathbf{i}\mathrm{Ric}\left(\nabla_{x_{1}}y_{2}+\nabla_{y_{1}}x_{2}\right), \quad \forall  x_{1},x_{2}, y_{1},y_{2}\in \mathfrak{g}_{0}.
\end{eqnarray*}
This asserts that
$\widetilde{\nabla}_{x_{1}+\mathbf{i}y_{1}}\circ\widetilde{\mathrm{Ric}}=\widetilde{\mathrm{Ric}}\circ\widetilde{\nabla}_{x_{1}+\mathbf{i}y_{1}}$
and hence $(\mathfrak{g},\langle\cdot,\cdot\rangle')$ is Ricci-parallel.
\end{proof}
\begin{rem}
We mention that the method in Proposition \ref{3-prop-complex} can also be used to construct flat metric Lie algebras \cite{am2003,boul2019,milnor76} and pseudo-Riemannian algebraic Ricci solitons \cite{cr2022JGA,cr2022JGP,yan21,yd21} from old ones.
In other words, if $(\mathfrak{g}_{0},\langle\cdot,\cdot\rangle_{0})$ is a flat metric Lie algebra or an algebraic Ricci soliton, then the metric Lie algebra $(\mathfrak{g},\langle\cdot,\cdot\rangle')$ defined as in Proposition \ref{3-prop-complex} is also flat or an algebraic Ricci soliton.
\end{rem}
Finally in this section, by using Proposition \ref{3-prop-complex}, we present certain examples of left invariant Ricci-parallel metrics on Lie groups.
\begin{example}
Let $\mathfrak{g}_{0}$ be a simple Lie algebra. According to formula \eqref{2-formula-ric}, its Killing form $\mathbf{K}$ is Einstein with Einstein constant $-\frac{1}{4}$. For instance, $\mathfrak{g}_{0}=\mathfrak{sl}(n,\mathbb{R})$ the special linear Lie algebra,
$$\mathbf{K}(x,y)=2n\mathrm{tr}(xy),\quad \forall x,y\in \mathfrak{sl}(n,\mathbb{R}).$$
So the metric $\langle\cdot,\cdot\rangle$ on $\mathfrak{g}=\mathfrak{sl}(n,\mathbb{C})$ defined by
\begin{eqnarray*}
  \langle x_{1}+\mathbf{i}y_{1},x_{2}+\mathbf{i}y_{2}\rangle
   &=&-\frac{n}{\lambda^{2}+\mu^{2}}[\lambda\mathrm{tr}(x_{1}x_{2})-\lambda\mathrm{tr}(y_{1}y_{2})
   +\mu\mathrm{tr}(x_{1}y_{2})+\mu\mathrm{tr}(y_{1}x_{2})],
\end{eqnarray*}
$x_{1},x_{2},y_{1},y_{2}\in \mathfrak{sl}(n,\mathbb{R})$, is Ricci-parallel but not Einstein, $\lambda,\mu\in \mathbb{R}$, $\mu\neq 0$.
\end{example}

\begin{example}
Suppose that $(\mathfrak{g}_{0},\langle\cdot, \cdot\rangle_{0})$ is an  Euclidean metric non-unimodular solvable Lie algebra and $\mathfrak{n}=[\mathfrak{g}_{0},\mathfrak{g}_{0}]$, $\mathfrak{a}=\mathfrak{n}^{\perp}$ the orthogonal complement of $\mathfrak{n}$ relative to $\langle\cdot, \cdot\rangle_{0}$.
According to the works of Heber \cite{heb98} and Lauret \cite{lau01,lau10},
 $(\mathfrak{g}_{0},\langle\cdot, \cdot\rangle_{0})$ is Einstein with Einstein constant $c<0$ if and only if the following conditions are satisfied.
\begin{itemize}
   \item [(a)]  $\{D_Y,T_Y:Y \in \mathfrak{a}\}$ is a family of commuting derivations of $\mathfrak{g}_{0}$ vanishing on $\mathfrak{a}$, where $D_Y (\ne 0)$, $T_Y$ denote the symmetric part and the skew-symmetric part of $\textnormal{ad}\,Y$ relative to $\langle \cdot, \cdot\rangle_{0},$ respectively.
   \item [(b)] $\langle Y, Y\rangle_{0}=-\frac{1}{c}\mathrm{tr}(D^2_Y),$ for any $Y \in \mathfrak{a}$.
   \item [(c)] $(\mathfrak{n},\langle \cdot, \cdot\rangle_{\mathfrak{n}})$ is a Ricci nilsoliton with Ricci operator $\textnormal{Ric}_{\mathfrak{n}}=c\mathrm{Id}+D$,
        where $\langle \cdot, \cdot\rangle_{\mathfrak{n}}
       =\langle \cdot,\cdot\rangle_{0}|_{{\mathfrak{n}}\times{\mathfrak{n}}}$
        and $D=(D_Y)|_\mathfrak{n}$ for some $Y\in \mathfrak{a}$.
\end{itemize}
For instance, let $(\mathfrak{n}=H_{n},\langle\cdot,\cdot\rangle_{\mathfrak{n}})$ be a $(2n+1)$-dimensional Euclidean metric Heisenberg Lie algebra with an orthonormal basis $\{E_{i}, Z\}_{i=1}^{2n}$ such that
$$[E_{2i-1},E_{2i}]=\sqrt{\frac{2}{n+2}}Z,\quad i=1,2,\ldots,n.$$
By a direct computation, one obtains the Ricci operator  $\mathrm{Ric}_{\mathfrak{n}}$
of $(\mathfrak{n},\langle\cdot,\cdot\rangle_{\mathfrak{n}})$ with respect to the basis $\{E_{i},Z\}_{i=1}^{2n}$ is  as follows:
\begin{eqnarray*}
\mathrm{Ric}_{\mathfrak{n}}&=&\mathrm{diag}\left(-\frac{1}{n+2},-\frac{1}{n+2},\cdots,-\frac{1}{n+2},\frac{n}{n+2}\right)\\
&=&-\mathrm{Id}+\mathrm{diag}\left(\frac{n+1}{n+2},\frac{n+1}{n+2},\cdots,\frac{n+1}{n+2},\frac{2(n+1)}{n+2}\right).
\end{eqnarray*}
Clearly, $(\mathfrak{n},\langle\cdot,\cdot\rangle_{\mathfrak{n}})$ is a Ricci nilsoliton with $c=-1$.
Now let
$$D=\mathrm{diag}\left(\frac{n+1}{n+2},\frac{n+1}{n+2},\cdots,\frac{n+1}{n+2},\frac{2(n+1)}{n+2}\right)$$
and extend $(\mathfrak{n}=H_{n},\langle\cdot,\cdot\rangle_{\mathfrak{n}})$ to an  Einstein metric solvable Lie algebra
$(\mathfrak{g}_{0}=\mathbb{R}D+\mathfrak{n},\langle\cdot, \cdot\rangle_{0})$ as follows:
$$\langle D,D\rangle_{0}=\mathrm{tr}(D^2)=\frac{2(n+1)^{2}}{n+2}, \quad \langle D,\mathfrak{n}\rangle_{0}=0.$$

By Proposition \ref{3-prop-complex},  on the complex Lie algebra $\mathfrak{g}=\mathbb{C}D+\mathfrak{n}^{\mathbb{C}}$, the complexification of $\mathfrak{g}_{0}$,
the metric $\langle\cdot, \cdot\rangle$ defined by
\begin{gather*}
\left\langle\sum_{i=1}^{2n}(a_{i}+\mathbf{i}b_{i})E_{i}+(c_{1}+\mathbf{i}c_{2})Z+(d_{1}+\mathbf{i}d_{2})D,
\sum_{i=1}^{2n}(a'_{i}+\mathbf{i}b'_{i})E_{i}+(c'_{1}+\mathbf{i}c'_{2})Z+(d'_{1}+\mathbf{i}d'_{2})D\right\rangle\\
=\frac{-2}{\lambda^{2}+\mu^{2}}\Bigg[\sum_{i=1}^{2n}\left(\lambda a_{i}a'_{i}-\lambda b_{i}b'_{i}+\mu a_{i}b'_{i}
+\mu b_{i}a'_{i}\right)
+\lambda c_{1}c'_{1}-\lambda c_{2}c'_{2}+\mu c_{1}c'_{2}+\mu c_{2}c'_{1}\\
+\frac{2(n+1)^{2}}{n+2}(\lambda d_{1}d'_{1}-\lambda d_{2}d'_{2}+\mu d_{1}d'_{2}+\mu d_{2}d'_{1})\Bigg],
\end{gather*}
$a_{i}, a'_{i}, b_{i}, b'_{i}, c_{1}, c'_{1}, c_{2}, c'_{2}, d_{1}, d'_{1}, d_{2}, d'_{2}\in \mathbb{R}$, is Ricci-parallel but not Einstein, $\lambda,\mu\in \mathbb{R}$, $\mu\neq 0$.
\end{example}

\section{Ricci-parallel metric Lie algebras of type II}\label{sec4}
In this section, we study the Ricci-parallel metric Lie algebras of type II.
Keep the notation as in the first section. Let  $(\mathfrak{g}_{0},\langle\cdot,\cdot\rangle_{0})$  be a metric Lie algebra  and
$(\mathfrak{g},\langle \cdot,\cdot\rangle)$ be a double extension of $(\mathfrak{g}_{0},\langle\cdot,\cdot\rangle_{0})$ defined by
the following form:
$$\mathfrak{g}=\mathbb{R}u+\mathbb{R}v+\mathfrak{g}_{0},$$
 $$[u,e]=D(e)+\langle L,e\rangle_{0}v,\quad [e,e']=[e,e']_{0}+\langle K(e),e'\rangle_{0}v,
\quad [\mathfrak{g},v]=0,\quad e,e'\in \mathfrak{g}_{0},$$
$$\langle u,v\rangle=1,\quad \langle u,u\rangle=\langle v,v\rangle=0,$$
$$\langle u,\mathfrak{g}_{0}\rangle=\langle v,\mathfrak{g}_{0}\rangle=0, \quad
\langle \cdot,\cdot\rangle|_{\mathfrak{g}_{0}}=\langle\cdot,\cdot\rangle_{0},$$
where $D,K\in \mathrm{End}(\mathfrak{g}_{0})$, $L\in \mathfrak{g}_{0}$ is a vector.
Notice that  $K$ is skew-symmetric with respect to $\langle\cdot,\cdot\rangle_{0}$ and
the Jacobi identity of $\mathfrak{g}$ are equivalent to the following conditions:
\begin{itemize}
  \item $D$ is a derivation of $\mathfrak{g}_{0}$,
  \item $\langle L,[e,e']_{0}\rangle_{0}=\langle (KD+D^{*}K)(e),e'\rangle_{0}$ holds for all $e,e'\in \mathfrak{g}_{0}$.
\end{itemize}

To state the main result of this section, we need the following two lemmas.
\begin{lem}\label{4-lem-connection}
Notation as above. The Levi-Civita connection $\nabla$ of $(\mathfrak{g},\langle \cdot,\cdot\rangle)$ is given as follows.
\begin{gather*}
  \nabla_{x}v=\nabla_{v}x=0, \quad \nabla_{u}u=-L, \quad x\in \mathfrak{g}.\\
  \nabla_{u}e=\frac{1}{2}(D-D^{*}-K)(e)+\langle L,e\rangle_{0} v, \quad e\in \mathfrak{g}_{0}.\\
  \nabla_{e}u=-\frac{1}{2}(D+D^{*}+K)(e),\quad \quad e\in \mathfrak{g}_{0}.\\
  \nabla_{e}f=\widetilde{\nabla}_{e}f+\frac{1}{2}\langle (D+D^{*}+K)(e),f\rangle_{0} v,
  \quad e,f\in \mathfrak{g}_{0}.
\end{gather*}
Here $\widetilde{\nabla}$ denotes the  Levi-Civita connection  of $(\mathfrak{g}_{0},\langle \cdot,\cdot\rangle_{0})$.
\end{lem}
\begin{proof}
Notice that $[\mathfrak{g},v]=0$ and $\langle[\mathfrak{g},\mathfrak{g}],v\rangle=0$, we easily obtain that
$\nabla_{x}v=\nabla_{v}x=0$ holds for all $x\in\mathfrak{g}$.
Then for all $e\in\mathfrak{g}_{0}$, note that
\begin{eqnarray*}
\langle U(u,u),e\rangle=\langle [e,u],u\rangle=-\langle L,e\rangle_{0},
\end{eqnarray*}
we obtain $\nabla_{u}u=U(u,u)=-L$.

Moreover, since for all $e, f\in \mathfrak{g}_{0}$,
\begin{eqnarray*}
\langle U(u,e),u\rangle=\frac{1}{2}\langle [u,e],u\rangle=\frac{1}{2}\langle L,e\rangle_{0}
\end{eqnarray*}
and
\begin{eqnarray*}
\langle U(u,e),f\rangle&=&\frac{1}{2}\left( \langle [f,u],e\rangle+\langle [f,e],u\rangle  \right)\\
                           &=&\frac{1}{2}\left(-\langle D(f),e\rangle_{0}-\langle K(e),f\rangle_{0} \right)\\
                           &=&-\frac{1}{2}\left\langle (D^{*}+K)(e),f\right\rangle_{0},
\end{eqnarray*}
we have
\begin{eqnarray*}
 U(u,e)=\frac{1}{2}\langle L,e\rangle_{0} v-\frac{1}{2}(D^{\ast}+K)(e)
\end{eqnarray*}
and consequently
\begin{eqnarray*}
\nabla_{u}e&=&\frac{1}{2}[u,e]+U(u,e)\\
             &=&\frac{1}{2}(D-D^{*}-K)(e)+\langle L,e\rangle_{0} v.
\end{eqnarray*}
Hence
\begin{eqnarray*}
\nabla_{e}u&=&\nabla_{u}e+[e,u]\\
             &=&\frac{1}{2}(D-D^{*}-K)(e)+\langle L,e\rangle_{0} v-D(e)-\langle L,e\rangle_{0} v\\
             &=&-\frac{1}{2}(D+D^{*}+K)(e).
\end{eqnarray*}

Finally, for all $e,f,g\in \mathfrak{g}_{0}$, $\langle U(e,f),v\rangle=0$,
\begin{eqnarray*}
\langle U(e,f),u\rangle&=&\frac{1}{2}\left( \langle [u,e],f\rangle+\langle [u,f],e\rangle  \right)\\
                           &=&\frac{1}{2}\left(\langle D(e),f\rangle_{0}+\langle D(f),e\rangle_{0} \right)\\
                           &=&\frac{1}{2}\left\langle (D+D^{*})(e),f\right\rangle_{0},
\end{eqnarray*}
and
\begin{eqnarray*}
\langle U(e,f),g\rangle&=&\frac{1}{2}\left( \langle [g,e],f\rangle+\langle [g,f],e\rangle  \right)\\
&=&\frac{1}{2}\left( \langle [g,e]_{0},f\rangle_{0}+\langle [g,f]_{0},e\rangle _{0} \right).
\end{eqnarray*}
 We obtain
 \begin{eqnarray*}
\nabla_{e}f&=&\frac{1}{2}[e,f]+U(e,f)\\
&=&\widetilde{\nabla}_{e}f+\frac{1}{2}\langle K(e),f\rangle_{0}v
+\frac{1}{2}\left\langle (D+D^{*})(e),f\right\rangle_{0}v\\
&=& \widetilde{\nabla}_{e}f+\frac{1}{2}\langle (D+D^{*}+K)(e),f\rangle_{0} v.
\end{eqnarray*}
This completes the proof of the lemma.
\end{proof}
\begin{lem}\label{4-lem-ric}
Notation as above. The Ricci operator $\mathrm{Ric}$ of $(\mathfrak{g},\langle \cdot,\cdot\rangle)$ is given as follows.
\begin{gather*}
  \mathrm{Ric}(u)=\Delta+\Gamma v, \\
   \mathrm{Ric}(e)=\mathrm{Ric}_{0}(e)+\langle\Delta,e\rangle_{0}v,\quad e\in \mathfrak{g}_{0},\\
    \mathrm{Ric}(v)=0.
\end{gather*}
Here
$$\Gamma=-\frac{1}{2}\mathrm{tr}(D^{2})-\frac{1}{2}\mathrm{tr}(D^{*}D)
  -\frac{1}{4}\mathrm{tr}(K^{2})+\langle L,Z_{0}\rangle\in \mathbb{R}$$
and $\Delta\in \mathfrak{g}_{0}$ is determined by
\begin{equation*}
  \langle \Delta, e\rangle_{0}=-\frac{1}{2}\mathrm{tr}((D+D^{*})\circ \mathrm{ad}_{0}\,e)
  +\frac{1}{2}\langle (D-K)(Z_{0}),e\rangle_{0}
 -\frac{1}{4}\mathrm{tr}(KS^{0}_{e}),\quad \forall e\in \mathfrak{g}_{0}.
\end{equation*}
$Z_{0}$ denotes the mean curvature vector of $(\mathfrak{g}_{0},\langle\cdot,\cdot\rangle_{0})$
and $\mathrm{Ric}_{0}$ denotes the Ricci operator of $(\mathfrak{g}_{0},\langle\cdot,\cdot\rangle_{0})$.
\end{lem}
\begin{proof}
It follows directly from formula \eqref{2-formula-ric}, see also Proposition 4.2 of \cite{yd23}.
Let $Z_{0}$ be the mean curvature vector of $(\mathfrak{g}_{0},\langle\cdot,\cdot\rangle_{0})$,
 then  the mean curvature vector $Z_{\mathfrak{g}}$ of $(\mathfrak{g},\langle\cdot,\cdot\rangle)$ is equal to $Z_{0}+(\mathrm{tr}(D))v$.
 Notice that $S_{u}(\mathfrak{g})\subset \mathfrak{g}_{0}+\mathbb{R}v$, $S_{u}|_{\mathfrak{g}_{0}}=K$ and $S_{u}(v)=0$,  we have
 \begin{eqnarray*}
   \mathrm{ric}(u,u) &=& -\frac{1}{2}\mathbf{K}(u,u)-\langle [Z_{\mathfrak{g}},u],u\rangle-\frac{1}{2}\mathrm{tr}((\mathrm{ad}\,u)^{*}\mathrm{ad}\,u)
  +\frac{1}{4}\mathrm{tr}(S^{*}_{u}S_{u})\\
    &=&  -\frac{1}{2}\mathrm{tr}(D^{2})-\frac{1}{2}\mathrm{tr}(D^{*}D)
  -\frac{1}{4}\mathrm{tr}(K^{2})+\langle L,Z_{0}\rangle_{0}\\
  &=&\Gamma,
 \end{eqnarray*}
 and
\begin{eqnarray*}
  \mathrm{ric}(u,e) &=& -\frac{1}{2}\mathbf{K}(u,e)-\frac{1}{2}(\langle [Z_{\mathfrak{g}},u],e\rangle+\langle [Z_{\mathfrak{g}},e],u\rangle)-\frac{1}{2}\mathrm{tr}((\mathrm{ad}\,u)^{*}\mathrm{ad}\,e)
  +\frac{1}{4}\mathrm{tr}(S^{*}_{u}S_{e}) \\
  &=&-\frac{1}{2}\mathrm{tr}((D+D^{*})\circ \mathrm{ad}_{0}\,e)
  +\frac{1}{2}\langle (D-K)(Z_{0}),e\rangle_{0}
 -\frac{1}{4}\mathrm{tr}(KS^{0}_{e})\\
   &=&  \langle \Delta, e\rangle_{0}, \quad \forall e\in \mathfrak{g}_{0}.
\end{eqnarray*}
Moreover, $\mathrm{ric}(e,f)=\mathrm{ric}_{0}(e,f)$, $e,f\in \mathfrak{g}_{0}$,  and $\mathrm{ric}(v,\mathfrak{g})=0$. This completes the proof of the lemma.
\end{proof}
We obtain a necessary and sufficient condition for $(\mathfrak{g},\langle \cdot,\cdot\rangle)$ to be Ricci-parallel.
\begin{thm}
$(\mathfrak{g},\langle \cdot,\cdot\rangle)$ is Ricci-parallel if and only if
$(\mathfrak{g}_{0},\langle\cdot,\cdot\rangle_{0})$ is Ricci-parallel and the following conditions hold:
\begin{description}
  \item[$(C_{1})$] $\langle L,\Delta\rangle_{0}=0;$
  \item[$(C_{2})$] $\mathrm{Ric}_{0}(L)=-\frac{1}{2}(D-D^{*}-K)(\Delta);$
  \item[$(C_{3})$] $\mathrm{Ric}_{0}\circ(D-D^{*}-K)=(D-D^{*}-K)\circ\mathrm{Ric}_{0};$
  \item[$(C_{4})$] $(D+D^{*}-K)(\Delta)=0;$
  \item[$(C_{5})$] $\widetilde{\nabla}_{e}\Delta=-\frac{1}{2}\mathrm{Ric}_{0}\left((D+D^{*}+K)(e)\right), \quad \forall e\in \mathfrak{g}_{0}.$
\end{description}
\end{thm}
\begin{proof}
By Proposition \ref{2-prop-parallel} and Lemma \ref{4-lem-connection},
 $(\mathfrak{g},\langle \cdot,\cdot\rangle)$ is Ricci-parallel if and only if the following two equations hold:
\begin{eqnarray}
  \mathrm{Ric}\circ \nabla_{u} &=&  \nabla_{u}\circ\mathrm{Ric}, \label{paraequ1}\\
   \mathrm{Ric}\circ \nabla_{e} &=&  \nabla_{e}\circ\mathrm{Ric}, \quad \forall e\in \mathfrak{g}_{0}. \label{paraequ2}
\end{eqnarray}
Clearly, \eqref{paraequ1} is equivalent to the following equations
\begin{eqnarray}
  \mathrm{Ric}(\nabla_{u}u) &=&  \nabla_{u}(\mathrm{Ric}(u)), \label{paraequ1-1}\\
   \mathrm{Ric}(\nabla_{u}e) &=&  \nabla_{u}(\mathrm{Ric}(e)), \quad \forall e\in \mathfrak{g}_{0}. \label{paraequ1-2}
\end{eqnarray}
By Lemmas \ref{4-lem-connection} and \ref{4-lem-ric}, \eqref{paraequ1-1} becomes
\begin{equation*}
 \mathrm{Ric}_{0}(-L)+\langle -L,\Delta\rangle_{0}v =\mathrm{Ric}(-L)=\nabla_{u}(\Delta)
 =\frac{1}{2}(D-D^{*}-K)(\Delta)+\langle L,\Delta\rangle_{0}v,
\end{equation*}
which is equivalent to $(C_{1})$ and $(C_{2})$.

Also, \eqref{paraequ1-2} becomes
\begin{equation*}
\mathrm{Ric}\left(\frac{1}{2}(D-D^{*}-K)(e)\right)=\nabla_{u}(\mathrm{Ric}_{0}(e)+\langle \Delta,e\rangle_{0}v),
\end{equation*}
and hence
\begin{eqnarray*}
&&\frac{1}{2}\mathrm{Ric}_{0}\left((D-D^{*}-K)(e)\right)+\frac{1}{2}\langle(D-D^{*}-K)(e),\Delta\rangle_{0}v\\
&=&\frac{1}{2}(D-D^{*}-K)(\mathrm{Ric}_{0}(e))+\langle L,\mathrm{Ric}_{0}(e)\rangle_{0} v,
\end{eqnarray*}
which is equivalent to $(C_{2})$ and $(C_{3})$.
Consequently, \eqref{paraequ1} is equivalent to $(C_{1})$, $(C_{2})$ and $(C_{3})$.

For \eqref{paraequ2}, it is equivalent to the following equations
\begin{eqnarray}
  \mathrm{Ric}(\nabla_{e}u) &=&  \nabla_{e}(\mathrm{Ric}(u)), \label{paraequ2-1}\\
   \mathrm{Ric}(\nabla_{e}f) &=&  \nabla_{e}(\mathrm{Ric}(f)), \quad \forall e,f\in \mathfrak{g}_{0}. \label{paraequ2-2}
\end{eqnarray}
By a direct calculation, \eqref{paraequ2-1} becomes
\begin{equation*}
\mathrm{Ric}\left(-\frac{1}{2}(D+D^{*}+K)(e)\right)=\nabla_{e}(\Delta)
=\widetilde{\nabla}_{e}\Delta+\frac{1}{2}\langle(D+D^{*}+K)(e),\Delta\rangle_{0}v,
\end{equation*}
and consequently
\begin{eqnarray*}
&&-\frac{1}{2}\mathrm{Ric}_{0}\left((D+D^{*}+K)(e)\right)
+\left\langle\Delta,-\frac{1}{2}(D+D^{*}+K)(e)\right\rangle_{0}v\\
&=&\widetilde{\nabla}_{e}\Delta+\frac{1}{2}\langle(D+D^{*}+K)(e),\Delta\rangle_{0}v,
\end{eqnarray*}
which is equivalent to $(C_{4})$ and $(C_{5})$.

Also, \eqref{paraequ2-2} becomes
\begin{equation*}
\mathrm{Ric}_{0}(\widetilde{\nabla}_{e}f)+\langle \Delta,\widetilde{\nabla}_{e}f\rangle_{0}v
=\widetilde{\nabla}_{e}(\mathrm{Ric}_{0}(f))
+\frac{1}{2}\langle(D+D^{*}+K)(e),\mathrm{Ric}_{0}(f)\rangle_{0}v,
\end{equation*}
which is equivalent to $(C_{5})$ and $(\mathfrak{g}_{0},\langle\cdot,\cdot\rangle_{0})$ is Ricci-parallel.
This completes the proof of the theorem.
\end{proof}

Theorem \ref{1-thm-main2} is a special case of the following corollary.
\begin{cor}\label{4-cor-delta}
Assume $(\mathfrak{g}_{0},\langle\cdot,\cdot\rangle_{0})$ is Ricci-parallel and $\Delta=0$,
 then $(\mathfrak{g},\langle \cdot,\cdot\rangle)$ is Ricci-parallel if and only if the following conditions are satisfied:
\begin{eqnarray*}
 \left\{
\begin{aligned}
&\mathrm{Ric}_{0}(L)=0,\\
&\mathrm{Ric}_{0}\circ(D-D^{*}-K)=(D-D^{*}-K)\circ\mathrm{Ric}_{0},\\
&\mathrm{Ric}_{0}\circ(D+D^{*}+K)=0.
\end{aligned}\right.
 \end{eqnarray*}
\end{cor}
\begin{rem}
By Lemma \ref{4-lem-ric}, it is easily seen that $(\mathfrak{g},\langle\cdot,\cdot\rangle)$ is Ricci flat if and only if $(\mathfrak{g}_{0},\langle\cdot,\cdot\rangle_{0})$ is Ricci flat and $\Delta=\Gamma=0$.  Conversely, Boucetta and Tibssirte in \cite{bt20} proved that every Lorentz Einstein metric  nilpotent Lie algebra with degenerate center must be a double extension of Euclidean metric Abelian Lie algebra.
\end{rem}
\begin{cor}
Let $(\mathfrak{g}_{0},\langle\cdot,\cdot\rangle_{0})$ be an ad-invariant metric nilpotent Lie algebra. Assume  $K=0$ and $L=0$, then $(\mathfrak{g},\langle \cdot,\cdot\rangle)$ is Ricci-parallel.
\end{cor}
\begin{proof}
Notice that $Z_{0}=0$, as $\mathfrak{g}_{0}$ is nilpotent.
Moreover, since $(\mathfrak{g}_{0},\langle\cdot,\cdot\rangle_{0})$ is ad-invariant, we have
$(\mathrm{ad}_{0}\,e)^{*}=-\mathrm{ad}_{0}\,e$ for all $e\in \mathfrak{g}_{0}$.
Hence
$$\mathrm{tr}(D^{*}\circ \mathrm{ad}_{0}\,e)=\mathrm{tr}(D^{*}\circ \mathrm{ad}_{0}\,e)^{*}
=-\mathrm{tr}(\mathrm{ad}_{0}\,e\circ D)=-\mathrm{tr}(D\circ \mathrm{ad}_{0}\,e),$$
and then $\langle \Delta, e\rangle_{0}=0$ holds for all $e\in \mathfrak{g}_{0}$.
It follows from Corollary \ref{4-cor-delta} and the fact that $\mathrm{Ric}_{0}=0$.
\end{proof}

Finally in this section, we give some interesting examples.
\begin{example}\label{4-exam-two-step}
Let $\mathfrak{g}_{0}$ be an Abelian Lie algebra endowed with an indefinite metric $\langle\cdot,\cdot\rangle_{0}$ of signature $(p,q)$.
Assume   $\mathfrak{D}\subset \mathrm{Der}(\mathfrak{g}_{0})$ is an $n$-dimensional Abelian Lie subalgebra
and  $\alpha\in Z^{2}(\mathfrak{D},\mathfrak{g}_{0})$ is a cocycle, that is,
 a skew-symmetric bilinear form from  $\mathfrak{D}\times \mathfrak{D}$ to $\mathfrak{g}_{0}$.
We can define a Lie bracket $[\cdot,\cdot]'$ on the vector space $\mathfrak{D}+\mathfrak{g}_{0}$  by
$$[D_{1}+e_{1}, D_{2}+e_{2}]'=D_{1}(e_{2})-D_{2}(e_{1})+\alpha(D_{1},D_{2}),$$
$D_{1},D_{2}\in \mathfrak{D}$, $e_{1},e_{2}\in \mathfrak{g}_{0}$,
where $\alpha$ should satisfy
$$D_{1}(\alpha(D_{2},D_{3}))+D_{2}(\alpha(D_{3},D_{1}))+D_{3}(\alpha(D_{1},D_{2}))=0,
\quad \forall D_{1},D_{2},D_{3}\in \mathfrak{D}.$$

Moreover, let $\mathfrak{D}^{*}$ be the dual space of $\mathfrak{D}$ and
$\theta\in Z^{2}((\mathfrak{D}+\mathfrak{g}_{0},[\cdot,\cdot]'), \mathfrak{D}^{*})$ be a cocycle.
Let  $\mathfrak{g}=(\mathfrak{D}+\mathfrak{g}_{0}+\mathfrak{D}^{*},[\cdot,\cdot])$
be the central extension of $(\mathfrak{D}+\mathfrak{g}_{0},[\cdot,\cdot]')$ by $\theta$:
$$[D_{1}+e_{1}+f_{1},D_{2}+e_{2}+f_{2}]=[D_{1}+e_{1},D_{2}+e_{2}]'+\theta(D_{1}+e_{1},D_{2}+e_{2}),$$
$D_{1},D_{2}\in \mathfrak{D}, e_{1},e_{2}\in \mathfrak{g}_{0}, f_{1},f_{2}\in \mathfrak{D}^{*}$.
We define the nature metric $\langle\cdot,\cdot\rangle$ of signature $(n+p,n+q)$ on $\mathfrak{g}$  as follows:
$$\langle D_{1}+e_{1}+f_{1},D_{2}+e_{2}+f_{2}\rangle=\langle e_{1},e_{2}\rangle_{0}
  +f_{1}(D_{2})++f_{2}(D_{1}),$$
$\forall D_{1},D_{2}\in \mathfrak{D}$, $e_{1},e_{2}\in \mathfrak{g}_{0}$,
  $f_{1},f_{2}\in \mathfrak{D}^{*}.$
  Notice that $\mathfrak{g}$ is a solvable Lie algebra,
  $[\mathfrak{g},\mathfrak{g}]\subset \mathfrak{g}_{0}+\mathfrak{D}^{*}$ and $\mathfrak{D}^{*}\subset C(\mathfrak{g})$, the center of $\mathfrak{g}$.
By definition,  we easily obtain that
 \begin{equation*}
   \nabla_{\mathfrak{g}}\mathfrak{g}\subset \mathfrak{g}_{0}+\mathfrak{D}^{*},\quad
   \nabla_{\mathfrak{g}}\mathfrak{D}^{*}=0,\quad
   \mathrm{Ric}(\mathfrak{g})\subset \mathfrak{D}^{*},\quad
   \mathrm{Ric}(\mathfrak{g}_{0}+\mathfrak{D}^{*})=0.
 \end{equation*}
 This implies that
 $\nabla_{\mathfrak{g}}\circ \mathrm{Ric}=\mathrm{Ric}\circ\nabla_{\mathfrak{g}}=0$.
 Consequently, $(\mathfrak{g},\langle\cdot,\cdot\rangle)$ is Ricci-parallel.
\end{example}

Note that every two step nilpotent Lie algebra can be expressed as the form of Example \ref{4-exam-two-step}, so we have
\begin{thm}
Every connected and simply connected two step nilpotent Lie group admits a left invariant pseudo-Riemannian Ricci-parallel metric.
\end{thm}
\begin{rem}
Up to know, it is still unknown whether every connected and simply connected two step nilpotent Lie group admits a left invariant pseudo-Riemannian Einstein metric, see \cite{yd23}.
\end{rem}

\begin{example}
We start with a class of  ad-invariant metric  Lie algebras,
 which is due to Bordemann \cite{bordemann97}.
Let $\mathfrak{D}$ be a Lie algebra with dual space $\mathfrak{D}^{*}$ and  $\theta\in Z^{2}(\mathfrak{D},\mathfrak{D}^{*})$ be a cocycle satisfying
$\theta(x_{1},x_{2})(x_{3})+\theta(x_{1},x_{3})(x_{2})=0$ for all $x_{1}, x_{2}, x_{3}\in \mathfrak{D}$.
Define  a Lie algebra $\mathfrak{g}_{0}=\mathfrak{D}+\mathfrak{D}^{*}$, direct sum as a vector space, where the Lie bracket $[\cdot,\cdot]_{0}$ is given by
$$[x_{1}+f_{1},x_{2}+f_{2}]_{0}=[x_{1},x_{2}]_{\mathfrak{D}}+x_{1}\cdot f_{2}-x_{2}\cdot f_{1}+\theta(x_{1},x_{2}),  \quad x_{1}, x_{2}\in \mathfrak{D}, f_{1},f_{2}\in \mathfrak{D}^{*},$$
where $x\cdot f$ denotes the coadjoint representation $\mathrm{ad}^{*}: \mathfrak{D}\rightarrow \mathrm{End}(\mathfrak{D}^{*})$,
 namely, $(x\cdot f)(y)=\mathrm{ad}^{*}(x)(f)(y)=-f([x,y]_{\mathfrak{D}})$, $ x,y\in \mathfrak{D}$, $f\in \mathfrak{D}^{*}$.
Then the natural  metric $\langle\cdot,\cdot\rangle_{0}$ on $\mathfrak{g}_{0}=\mathfrak{D}+\mathfrak{D}^{*}$
$$\langle x_{1}+f_{1}, x_{2}+f_{2}\rangle_{0}=f_{1}(x_{2})+f_{2}(x_{1})$$
is ad-invariant. In particular, $\mathfrak{g}_{0}$ is nilpotent when $\mathfrak{D}$ is nilpotent.
Conversely, Bordemann \cite{bordemann97} proved that every even-dimensional nilpotent Lie algebra  with an ad-invariant metric
is of this form.

To construct Ricci-parallel metric Lie algebras from the above ad-invariant metric  Lie algebras $(\mathfrak{g}_{0},\langle\cdot,\cdot\rangle_{0})$, we now  assume $\mathfrak{D}$ is an Abelian Lie algebra  endowed with an indefinite metric $\langle\cdot,\cdot\rangle_{\mathfrak{D}}$ and in this case
$\mathfrak{g}_{0}$ is at most two step nilpotent.
Let $\varphi:\mathfrak{D}\rightarrow \mathfrak{D}^{*}$ be the linear isomorphism given by $\varphi(x)(y)=\langle x,y\rangle_{\mathfrak{D}}$, $\forall x,y\in \mathfrak{D}$.
Moreover, let $D=0$, $L=0$ and $K$ be the skew-symmetric endomorphism on $(\mathfrak{g}_{0},\langle\cdot,\cdot\rangle_{0})$
given by
$$K(x)=K_{\mathfrak{D}}(x),\quad K(f)=\varphi(K_{\mathfrak{D}}(\varphi^{-1}(f))),\quad \forall x\in\mathfrak{D}, f\in \mathfrak{D}^{*},$$
where $K_{\mathfrak{D}}$ denotes a skew-symmetric endomorphism on $\mathfrak{D}$ with respect to $\langle\cdot,\cdot\rangle_{\mathfrak{D}}$.
It is easily seen that
$$ \langle \Delta, e\rangle_{0}=-\frac{1}{4}\mathrm{tr}(KS^{0}_{e})
 =-\frac{1}{4}\mathrm{tr}(K\circ \mathrm{ad}_{0}\,e)=0,\quad \forall e\in \mathfrak{g}_{0},$$
 since $(K\circ \mathrm{ad}_{0}\,e)(\mathfrak{D}) \subset \mathfrak{D}^{*}$
 and $(K\circ \mathrm{ad}_{0}\,e)(\mathfrak{D}^{*})=0$.
 This shows that the metric Lie algebra $(\mathfrak{g},\langle\cdot,\cdot\rangle)$ obtained by the double extension of $(\mathfrak{g}_{0},\langle\cdot,\cdot\rangle_{0})$ with respect to $(D=0,L=0,K)$ is Ricci-parallel.
 In particular, according to Lemma \ref{4-lem-ric},  $(\mathfrak{g},\langle\cdot,\cdot\rangle)$ is Ricci flat if and only if $\mathrm{tr}(K_{\mathfrak{D}}^{2})=0$.
\end{example}

\section{Lorentz Ricci-parallel  metric Lie algebras}\label{sec5}

In this section, we study the structure of  Lorentz Ricci-parallel  metric Lie algebras of type II.
Let  $(\mathfrak{g},\langle\cdot,\cdot\rangle)$ be an $(n+2)$-dimensional Lorentz Ricci-parallel  metric Lie algebra such that $\mathrm{Ric}^{2}=0$ but $\mathrm{Ric}\neq 0$.  Then according to a result of O'Neill,
 see Exercise 9.19 of \cite{oneill1983book}, there exists a basis $\{u,v,e_{1},\ldots,e_{n}\}$ with all scalar products zero except
$\langle u,v\rangle=\langle e_{1},e_{1}\rangle=\cdots=\langle e_{n},e_{n}\rangle=1$,
such that the matrix of $\mathrm{Ric}$ with respect to this basis is of the form
\begin{eqnarray*}
\left(\begin{array}{cc}
0 & 0\\
1 & 0
\end{array}\right)\oplus O_{n}.
\end{eqnarray*}
It is known that the matrix of a skew-symmetric endomorphism of $(\mathfrak{g},\langle\cdot,\cdot\rangle)$
with respect to the above basis is of the form
\begin{eqnarray*}
\left(\begin{array}{ccc}
a       &  0   &   -\beta^{\mathrm{T}}\\
0       & -a   &   -\alpha^{\mathrm{T}}\\
\alpha  &\beta &     A
\end{array}\right),
\end{eqnarray*}
where $a\in \mathbb{R}$, $\alpha,\beta\in \mathbb{R}^{n}$ are column vectors, $A$ is an $n\times n$ skew-symmetric matrix.
Now by Proposition \ref{2-prop-parallel}, for every $x\in \mathfrak{g}$,  the equality
$\nabla_{x}\circ \mathrm{Ric}=\mathrm{Ric}\circ \nabla_{x}$
implies that the matrix of $\nabla_{x}$ is of the form
\begin{eqnarray*}
\left(\begin{array}{ccc}
0       &   0 &     0\\
0       &   0 &   -\alpha^{\mathrm{T}}\\
\alpha  &   0 &       A
\end{array}\right), \quad \alpha \in \mathbb{R}^{n}.
\end{eqnarray*}
Hence we have
\begin{lem}
$[\mathfrak{g},\mathfrak{g}]\subset \mathbb{R}v+\mathfrak{g}_{0}$. Consequently, $\mathbb{R}v+\mathfrak{g}_{0}$ is an ideal of $\mathfrak{g}$, where $\mathfrak{g}_{0}=\mathbb{R}e_{1}+\cdots+\mathbb{R}e_{n}$ is the subspace of $\mathfrak{g}$.
\end{lem}
\begin{proof}
For all $x,y\in \mathfrak{g}$, we have $[x,y]=\nabla_{x}y-\nabla_{y}x\in \mathbb{R}v+\mathfrak{g}_{0}$.
\end{proof}
Moreover, notice that $\nabla_{x}v=0$ holds for all $x\in\mathfrak{g}$
and $\langle v, [\mathfrak{g},\mathfrak{g}]\rangle=0$. We obtain the following
\begin{lem}
For all $ x,y\in\mathfrak{g}$, one has
\begin{eqnarray*}
\langle [v,x],y\rangle+\langle [v,y],x\rangle=0.
\end{eqnarray*}
\end{lem}
\begin{proof}
Since $0=\nabla_{x}v=\frac{1}{2}[x,v]+U(x,v)$, we have
\begin{eqnarray*}
  0 &=& \langle\frac{1}{2}[x,v]+U(x,v),y\rangle \\
   &=& \frac{1}{2}\left(\langle [x,v],y\rangle+\langle [y,v],x\rangle\right).
\end{eqnarray*}
This completes the proof.
\end{proof}

\textbf{Proof of Theorem \ref{1-thm-main3}}.
By assumptions and formula \eqref{2-formula-ric}, one  has
\begin{eqnarray*}
0=\mathrm{ric}(v,v)=-\frac{1}{2}\sum_{i=1}^{n}\langle[v,e_{i}],[v,e_{i}]\rangle,
\end{eqnarray*}
we obtain that $[v,e_{i}]\in \mathbb{R}v$, $i=1,\ldots,n$.
Since $\mathfrak{g}$ is a nilpotent Lie algebra, we have $[v,e_{i}]=0$ for all $e_{i}$ and consequently $[v,\mathbb{R}v+\mathfrak{g}_{0}]=0$.
Therefore there exist a unique Lie bracket $[\cdot,\cdot]_{0}$ on $\mathfrak{g}_{0}$ and a unique cocycle $\theta\in Z^{2}(\mathfrak{g}_{0},\mathbb{R})$ such that
$$[e,f]=[e,f]_{0}+\theta(e,f)v, \quad \forall e,f\in \mathfrak{g}_{0}.$$
Clearly, $(\mathfrak{g}_{0},[\cdot,\cdot]_{0})$ is isomorphic to the quotient Lie algebra $(\mathbb{R}v+\mathfrak{g}_{0},[\cdot,\cdot])/\mathbb{R}v$.
Also, note that
\begin{eqnarray*}
0=\mathrm{ric}(u,v)&=&
-\frac{1}{2}\left(\left\langle\left[u,\frac{u+v}{\sqrt{2}}\right],\left[v,\frac{u+v}{\sqrt{2}}\right]\right\rangle
-\left\langle\left[u,\frac{u-v}{\sqrt{2}}\right],\left[v,\frac{u-v}{\sqrt{2}}\right]\right\rangle\right)\\
                   &=&\frac{1}{2}\langle [u,v],[u,v]\rangle,
\end{eqnarray*}
we have $[u,v]\in \mathbb{R}v$ and consequently $[u,v]=0$, as $\mathfrak{g}$ is nilpotent. So $[v,\mathfrak{g}]=0$ and $(\mathfrak{g}, \langle\cdot,\cdot\rangle)$ is a double extension of the Euclidean metric Lie algebra $(\mathfrak{g}_{0},[\cdot,\cdot]_{0},\langle\cdot,\cdot\rangle_{0})$, where $\langle\cdot,\cdot\rangle_{0}=\langle\cdot,\cdot\rangle|_{\mathfrak{g}_{0}}$.

 Finally, we show that $(\mathfrak{g}_{0},[\cdot,\cdot]_{0})$ is  Abelian.
 Indeed, by Lemma \ref{4-lem-ric}, the Ricci curvature $\mathrm{ric}_{0}$ of $(\mathfrak{g}_{0},[\cdot,\cdot]_{0},\langle\cdot,\cdot\rangle_{0})$ is equal to zero, which implies that $\mathfrak{g}_{0}$ is an Abelian Lie algebra, according to a result of Milnor \cite{milnor76}.
This completes the proof of  Theorem \ref{1-thm-main3}.\qed

At the last of this section, we will show that every Lorentz Ricci-parallel  metric Lie algebra of type II of dimensions 3 and 4 is a double extension of Euclidean metric Abelian Lie algebra.
\begin{thm}
 Every Lorentz Ricci-parallel  metric Lie algebra of type II of dimension 3  is a double extension of Euclidean metric Abelian Lie algebra.
\end{thm}
\begin{proof}
Since $\dim \mathfrak{g}=3$,
the matrix of $\mathrm{ad}\,v=\nabla_{v}$ with respect to the basis $\{u,v,e_{1}\}$ is of the form
 \begin{eqnarray*}
\left(\begin{array}{ccc}
0       &   0 &     0\\
0       &   0 &   -\gamma\\
\gamma  &   0 &       0
\end{array}\right)
\end{eqnarray*}
for some $\gamma\in \mathbb{R}$.
So we have $[v,e_{1}]=-\gamma v$ and $[v,u]=\gamma e_{1}$. Moreover,
\begin{eqnarray*}
  \gamma^{2}e_{1}=\gamma [v,u] &=& [u,[v,e_{1}]] \\
   &=& [[u,v],e_{1}]+[v,[u,e_{1}]] \\
   &=& [v,[u,e_{1}]]\in \mathbb{R}v,
\end{eqnarray*}
since $[u,e_{1}]\in \mathbb{R}v+\mathbb{R}e_{1}$.
We obtain that $\gamma=0$ and consequently $(\mathfrak{g}, \langle\cdot,\cdot\rangle)$ is a double extension of $(\mathbb{R}e_{1},\langle\cdot,\cdot\rangle|_{\mathbb{R}e_{1}})$.
\end{proof}
\begin{thm}
 Every Lorentz Ricci-parallel  metric Lie algebra of type II of dimension 4  is a double extension of Euclidean metric Abelian Lie algebra.
\end{thm}
\begin{proof}
We claim that $[v,\mathfrak{g}]=0$ and in this case
$(\mathfrak{g}, \langle\cdot,\cdot\rangle)$ is a double extension of $(\mathfrak{g}_{0},\langle\cdot,\cdot\rangle|_{\mathfrak{g}_{0}})$.
Since $\dim \mathfrak{g}_{0}=2$ and $(\mathfrak{g}_{0},\langle\cdot,\cdot\rangle|_{\mathfrak{g}_{0}})$ is Ricci flat, it is easy to see that $\mathfrak{g}_{0}$ is Abelian.

Now we suppose conversely that $\mathrm{ad}\,v\neq 0$.
Then there exist three constants $\gamma, \gamma_{1}, \gamma_{2}\in \mathbb{R}$ such that the matrix of $\mathrm{ad}\,v=\nabla_{v}$ with respect to the basis $\{u,v,e_{1},e_{2}\}$ is of the form
 \begin{eqnarray*}
\left(\begin{array}{cccc}
0       &   0 & 0  &  0\\
0       &   0 & -\gamma_{1}  &-\gamma_{2}\\
\gamma_{1}  &   0 &  0  &   -\gamma\\
\gamma_{2}   &   0 &  \gamma  & 0
\end{array}\right).
\end{eqnarray*}
We have two cases: $\gamma=0$, $\gamma_{1}^{2}+\gamma_{2}^{2}\neq 0$ and $\gamma\neq0$.

\textbf{Case 1}. $\gamma=0$, $\gamma_{1}^{2}+\gamma_{2}^{2}\neq 0$.
Without loss of generality we can assume $\gamma_{1}\neq 0$.
Then let $e'_{2}=e_{2}-\frac{\gamma_{2}}{\gamma_{1}}e_{1}$, we obtain
$$[v,e'_{2}]=\left[v,e_{2}-\frac{\gamma_{2}}{\gamma_{1}}e_{1}\right]
=[v,e_{2}]-\frac{\gamma_{2}}{\gamma_{1}}[v,e_{1}]=0.$$
Since $[v,[e_{1},e'_{2}]]=[[v,e_{1}],e'_{2}]+[e_{1},[v,e'_{2}]]=0$, we have
$[e_{1},e'_{2}]\in \mathbb{R}v+\mathbb{R}e'_{2}$. Therefore $\mathbb{R}v+\mathbb{R}e'_{2}$ is an Abelian ideal of $\mathbb{R}v+\mathfrak{g}_{0}$. In particular, $\mathbb{R}v+\mathfrak{g}_{0}$ is a non-nilpotent solvable Lie algebra, because $[e_{1},v]=\gamma_{1}v$.
 Now a direct calculation shows the following
 \begin{eqnarray*}
   [u,v]=-[v,u] &=& -\gamma_{1}e_{1}-\gamma_{2}e_{2} \\
    &=& -\gamma_{1}e_{1}-\gamma_{2}\left(e'_{2}+\frac{\gamma_{2}}{\gamma_{1}}e_{1}\right) \\
    &=& -\gamma_{2}e'_{2}-\frac{\gamma_{1}^{2}+\gamma_{2}^{2}}{\gamma_{1}}e_{1} \notin \mathbb{R}v+\mathbb{R}e'_{2},
 \end{eqnarray*}
which is in contradiction with the fact that every derivation of $\mathbb{R}v+\mathfrak{g}_{0}$ maps $\mathbb{R}v+\mathfrak{g}_{0}$ to
$\mathbb{R}v+\mathbb{R}e'_{2}$,
see Theorem 3.7 of \cite{Jaco1962}.

\textbf{Case 2}. $\gamma\neq0$.
In this case, let $e'_{1}=e_{1}+\frac{\gamma_{2}}{\gamma}v$ and $e'_{2}=e_{2}-\frac{\gamma_{1}}{\gamma}v$,
then we have
\begin{equation*}
  [v,e'_{1}]=[v,e_{1}]=-\gamma_{1}v+\gamma e_{2}=\gamma e'_{2}
\end{equation*}
and
\begin{equation*}
  [v,e'_{2}]=[v,e_{2}]=-\gamma_{2}v-\gamma e_{1}=-\gamma e'_{1}.
\end{equation*}
Moreover, let
$$u'=u-\frac{\gamma_{1}^{2}+\gamma_{2}^{2}}{2\gamma^{2}}v
-\frac{\gamma_{2}}{\gamma}e'_{1}+\frac{\gamma_{1}}{\gamma}e'_{2},$$
we obtain that
\begin{equation*}
  \langle u',u'\rangle=0,\quad \langle u',v\rangle=1,\quad \langle u',e'_{1}\rangle=\langle u',e'_{2}\rangle=0.
\end{equation*}
Hence with respect to the new basis $\{u',v,e'_{1},e'_{2}\}$, the matrices of $\mathrm{Ric}$ and $\mathrm{ad}\,v=\nabla_{v}$ are of the forms
  \begin{eqnarray*}
\left(\begin{array}{cccc}
0       &   0 & 0  &  0\\
1       &   0 & 0  & 0\\
0       &   0 &  0  & 0\\
0       &   0 &  0  & 0
\end{array}\right)
\end{eqnarray*}
and
 \begin{eqnarray*}
\left(\begin{array}{cccc}
0       &   0 & 0  &  0\\
0       &   0 & 0  & 0\\
0       &   0 &  0  &   -\gamma\\
0        &   0 &  \gamma  & 0
\end{array}\right),
\end{eqnarray*}
respectively. In particular, $[v,u']=0$.
Now notice that
\begin{equation*}
  [v,[e'_{1},e'_{2}]]=[[v,e'_{1}],e'_{2}]+[e'_{1},[v,e'_{2}]]=0,
\end{equation*}
one has  $[e'_{1},e'_{2}]\in \mathbb{R}v$.

If $[e'_{1},e'_{2}]=0$, the Lie algebra $\mathbb{R}v+\mathfrak{g}_{0}=\mathbb{R}v+\mathbb{R}e'_{1}+\mathbb{R}e'_{2}$ is solvable. Moreover, $\mathrm{ad}\,u'$ maps $\mathbb{R}v+\mathfrak{g}_{0}$ to $\mathbb{R}e'_{1}+\mathbb{R}e'_{2}$, the nilradical of $\mathbb{R}v+\mathfrak{g}_{0}$.
Hence $\mathrm{ad}\,u'|_{\mathbb{R}v+\mathfrak{g}_{0}}=\beta (\mathrm{ad}\,v|_{\mathbb{R}v+\mathfrak{g}_{0}})$
for some constant $\beta\in \mathbb{R}$, which implies that $\mathfrak{g}$ is unimodular.
Note that
$\langle u',[\mathfrak{g},\mathfrak{g}]\rangle=0$,
we obtain
\begin{eqnarray*}
  1 &=& \mathrm{ric}(u',u') \\
   &=& -\frac{1}{2}\mathbf{K}(u',u')-\frac{1}{2}\langle[u',e'_{1}],[u',e'_{1}]\rangle
   -\frac{1}{2}\langle[u',e'_{2}],[u',e'_{2}]\rangle \\
   &=& -\frac{1}{2}\beta^{2}(-2\gamma^{2})-\frac{1}{2}\beta^{2}\gamma^{2}-\frac{1}{2}\beta^{2}\gamma^{2} \\
   &=& 0,
\end{eqnarray*}
which is a contradiction.

On the other hand, if $[e'_{1},e'_{2}]=\xi v$ for some $\xi\neq0$,
then the Lie algebra $\mathbb{R}v+\mathfrak{g}_{0}=\mathbb{R}v+\mathbb{R}e'_{1}+\mathbb{R}e'_{2}$ is simple and  thus every derivation is inner. Hence  $\mathrm{ad}\,u'|_{\mathbb{R}v+\mathfrak{g}_{0}}=\beta (\mathrm{ad}\,v|_{\mathbb{R}v+\mathfrak{g}_{0}})$
for some constant $\beta\in \mathbb{R}$, which also implies that $\mathfrak{g}$ is unimodular.
Now we have
\begin{eqnarray*}
  0 &=& \mathrm{ric}(e'_{1},e'_{1}) \\
   &=& -\frac{1}{2}\mathbf{K}(e'_{1},e'_{1})
   -\frac{1}{2}\left(\left\langle\left[e'_{1},\frac{u'+v}{\sqrt{2}}\right],
   \left[e'_{1},\frac{u'+v}{\sqrt{2}}\right]\right\rangle
   -\left\langle\left[e'_{1},\frac{u'-v}{\sqrt{2}}\right],
   \left[e'_{1},\frac{u'-v}{\sqrt{2}}\right]\right\rangle\right) \\
   &&+\frac{1}{2}\left(\left\langle\left[e'_{2},\frac{u'+v}{\sqrt{2}}\right],e'_{1}\right\rangle^{2}
   -\left\langle\left[e'_{2},\frac{u'-v}{\sqrt{2}}\right],e'_{1}\right\rangle^{2}\right)  \\
   &=& -\frac{1}{2}(-2\gamma\xi)
   -\frac{1}{2}\left(\frac{1}{2}(\beta+1)^{2}\gamma^{2}-\frac{1}{2}(\beta-1)^{2}\gamma^{2}\right)
   +\frac{1}{2}\left(\frac{1}{2}(\beta+1)^{2}\gamma^{2}-\frac{1}{2}(\beta-1)^{2}\gamma^{2}\right)\\
   &=&\gamma\xi,
\end{eqnarray*}
which is also a contradiction. This completes the proof of the theorem.
\end{proof}


\end{document}